\documentclass[a4paper,12pt]{amsart}
\usepackage[dvipdfm]{color}
\usepackage[dvipdfm,bookmarks=true,bookmarksnumbered=true,bookmarkstype=toc]{hyperref}

\usepackage{amssymb}
\usepackage{amsxtra}
\numberwithin{equation}{section}

\newtheorem{thm}{Theorem}[section]
\newtheorem{lem}[thm]{Lemma}
\newtheorem{cor}[thm]{Corollary}
\newtheorem{prop}[thm]{Proposition}

\theoremstyle{definition}
\newtheorem{defn}[thm]{Definition}

\theoremstyle{remark}
\newtheorem{rem}[thm]{Remark}

\newcommand{\GL}{\operatorname{GL}}
\newcommand{\Ind}{\operatorname{Ind}}
\newcommand{\ind}{\operatorname{ind}}

\newcommand{\Tr}{\operatorname{Tr}}

\newcommand{\of}{{\mathcal O}_{F}}
\newcommand{\Oe}{{\mathcal O}_{E}}

\newcommand{\calL}{\mathcal L}
\newcommand{\Ker}{\operatorname{Ker}}
\newcommand{\tr}{\operatorname{tr}}
\newcommand{\n}{\operatorname{N}}
\newcommand{\M}{\operatorname{M}}
\newcommand{\End}{\operatorname{End}}
\newcommand{\Aut}{\operatorname{Aut}}
\newcommand{\sgn}{\operatorname{sgn}}
\newcommand{\ch}{\operatorname{ch}}

\newcommand{\fA}{\mathfrak A}
\newcommand{\fP}{\mathfrak P}
\newcommand{\fK}{\mathfrak K}

\begin{document}
\title[Epsilon factor of $\GL_l\times \GL_{l'}$]{Epsilon factor for $\GL_{l}\times \GL_{l'}$; $l\neq l'$ primes}
\author{Tetsuya Takahashi}
\address{Department of Mathematics and Information Science, College of Integrated Arts and Sciences, Osaka Prefecture University}
\email{takahasi@mi.cias.osakafu-u.ac.jp}
\keywords{$\epsilon$-factor of pairs, supercuspidal representation, local Langlands correspondence} 
\subjclass{Primary 22E50, Secondary 11F70}
\begin{abstract}
Let $F$ be a non-Archimedean local field with finite residual field of characteristic $p$.
In this article we calculate the $\varepsilon$-factor of pairs for 
$\GL_l(F) \times \GL_{l'}(F)$ where $l$ and $l'$ are distinct primes including the case $l=p$.
For this calculation, we use the local Langlands correspondence and non-Galois base change lift.
This method leads to the explicit conjecture of the $\varepsilon$-factor of 
the representations of $\GL_m \times \GL_n$ when $n$ is relatively prime to $m$ and $p$.
\end{abstract}
\maketitle

\section{Introduction}
Let $F$ be a non-Archimedean local field with finite residual field of characteristic $p$ and
the ${\mathcal W}_F$ the absolute Weil group of $F$. For an integer $n \geq 1$, 
we denote by ${\mathcal A}_n(F)$ the set of equivalent classes of irreducible supercuspidal
representations of $\GL_n(F)$ and by ${\mathcal G}_n(F)$ the set of equivalent classes of irreducible continuous complex representations of ${\mathcal W}_F$ of dimension $n$.
The local Langlands conjecture tells us that there exists a unique bijection $\Lambda_n^F$ from 
${\mathcal G}_n(F)$ to ${\mathcal A}_n(F)$ which satisfies the following conditions:
\begin{enumerate}
\item 
 For $\chi \in \widehat{F^{\times}}$and $\sigma\in{\mathcal G}_n(F)$, 
  \begin{equation}\label{LLL1}
  \Lambda_n^F(\chi\sigma)=\chi\Lambda_n^F(\sigma)
  \end{equation}
  (By the reciprocity map of local class field theory, we identify $\widehat{F^{\times}}={\mathcal A}_1(F)$ with $W_F^{ab}={\mathcal G}_1(F)$. By this identification, $\Lambda_1$ is the
  identity map.)
\item 
  For $\sigma\in{\mathcal G}_n(F)$, 
 \begin{equation}\label{LLL2}
  \Lambda_n^F(\check{\sigma})=\Lambda_n^F(\sigma)\spcheck .
 \end{equation}
 \item
  Let $\omega_\pi$ denote the central quasi-character of $\pi\in{\mathcal A}_n(F)$.  
  For $\sigma\in{\mathcal G}_n(F)$.
 \begin{equation}\label{LLL3}
  \omega_{\Lambda_n^F(\sigma)}=\det \sigma.
 \end{equation}
 \item
  Let $\psi_F$ be a non-trivial character of $F$. 
  For $\sigma\in{\mathcal G}_n(F)$, 
 \begin{equation}\label{LLL4}
  \varepsilon(\Lambda_n^F(\sigma),s,\psi_F)=\varepsilon(\sigma,s,\psi_F).
 \end{equation}
 where the left hand side is the Godement-Jacquet local constant \cite{GJ} and 
 the right hand side is the Langlands-Deligne local constant \cite{De}. (In fact, this condition
 is contained in the following condition (5).)
\item
For $\sigma\in{\mathcal G}_n(F)$ and $\sigma'\in{\mathcal G}_{n'}(F)$,
\begin{equation}\label{LLL5}
 \varepsilon(\Lambda_n^F(\sigma)\times \Lambda_{n'}^F(\sigma'),s,\psi_F)
 =\varepsilon(\sigma\otimes\sigma',s,\psi_F)
\end{equation}
where the $\varepsilon$-factor of pairs of the left hand side is in the sense of \cite{JPSS}, 
\cite{Sh}.
\end{enumerate}
This conjecture has been proved in \cite{LRS} when $\ch F=p$ and in \cite{HT}, \cite{He2} when
$\ch F=0$. But their proof relies on the global tool and say nothing explicit
about the local Langlands correspondence.

On the other hand, there are some explicit correspondences in the following cases:
\begin{enumerate}
\item
When $(n,p)=1$, Howe-Moy \cite{Ho},\cite{Mo} gives an explicit bijection between ${\mathcal G}_n(F)$ and ${\mathcal A}_n(F)$ when $(n,p)=1$. (See also \cite{Rei}).
\item
When $n=p$, Kutzko-Moy \cite{KM} gives an explicit bijection between ${\mathcal G}_n(F)$ and ${\mathcal A}_n(F)$ . (See also \cite{He1}).
\item
When $n=p^m$, Bushnell-Henniart \cite{BH2} gives an explicit bijection between ${\mathcal G}_{p^m}^{wr}(F)$ and ${\mathcal A}_{p^m}^{wr}(F)$ .
(For the definition of ${\mathcal G}_{p^m}^{wr}(F)$ and ${\mathcal A}_{p^m}^{wr}(F)$ , see below Remark~\ref{tamepi},)
\end{enumerate}
All three bijections satisfy the condition (1)--(4) of the local Langlands correspondence. Thus the
main obstacle to get the explicit local Langlands correspondence is $\varepsilon$-factor of pairs.
We know very few about the explicit calculation of $\varepsilon(\pi_1\times\pi_2,s,\psi_F)$ for
$\pi_1\in{\mathcal A}_{n_1}(F)$ and $\pi_2\in{\mathcal A}_{n_2}(F)$;
The known cases are (i) $n_1=n_2$ (\cite{Li}) and $\pi_2=\check{\pi}_1$ (\cite{BH4}),
(ii) $\pi_1\in {\mathcal A}_{p^{i_1}}^{wr}(F)$ and $\pi_2\in {\mathcal A}_{p^{i_2}}^{wr}(F)$
(\cite{BH5}).

In this paper we consider the case $n_1\neq n_2$ are primes. Set $n_1=l$ and $n_2=l'$.
We admit the case $l=p$. Since $l\neq l'$, we may assume $l'\neq p$.
We get the  relation of $\varepsilon$ -factor of $\GL_{l}(F)\times \GL_{l'}(F)$ with 
$\varepsilon$-factor of $\GL_l(E)$ where $E$ is an extension of $F$ associated with $\pi_2$.
(See Theorem~\ref{main}.) 
 
Let us summarize the contents of this paper, indicating its organization:

Section 1 reviews the construction of irreducible supercuspidal representations $\pi$
of $\GL_l(F)$ and the explicit formula of $\varepsilon(\pi,s,\psi_F)$. All of this section is well-known.  Section 2 is devoted to review some explicit correspondences and the tame lifting.
When $l\neq p$, ${\mathcal G}_l(F)$ consists of the representations in the form $\Ind_{W_E}^{W_F} \theta$; $E/F$ is an extension of degree $l$ and $\theta$ is a quasi-character of $E^{\times}$. By way of such $\theta$, there is very explicit Howe-Moy correspondence between 
${\mathcal G}_l(F)$ and ${\mathcal A}_l(F)$. But when $l=p$, there exists non-monomial 
representations in ${\mathcal G}_p(F)$; so we need the tame base change lift to get the
correspondence. (See \cite{KM},\cite{BH2}.)
Let $\pi\in{\mathcal A}_p^{wr}(F)$ and $K/F$ a tamely ramified extension. 
After the definition of \cite{BH},  we give the tame base change lift $l_{K}(\pi)$ explicitly (Theorem~\ref{explicit lift})
and show $l_{K}$ is compatible with the local Langlands correspondence
(Proposition~\ref{compatibility of lifting}).  We also define the tame base change lift 
$l_{K}$ for the case $l\neq p$ and prove it is compatible with the Howe-Moy correspondence
(Proposition~\ref{compatibility of lifting II}). These are essential tool to calculate 
the $\epsilon$-factor of pairs.
Section 3 calculates the $\varepsilon$-factor of $\GL_{l'}(F)\times\GL_{l}(F)$.
By the result of Bushnell-Henniart \cite{BH6},  
the Howe-Moy correspondence coincides with the Local Langlands correspondence for $\GL_{l}(F)$.
Thus we calculate the $\varepsilon$-factor in the Galois side and then transfer it to the automorphic
 side using the results in section 2.
 
 
\vspace{1cm}

{\bf\large Notation}

\vspace{8pt}
\par
  Let $F$ be a non-archimedean local field. 
  We denote by ${\mathcal O}_F$, $P_F$, $\varpi_F$, $k_F$ and $v_F$ 
 the maximal order of $F$, the maximal ideal of  ${\mathcal O}_F$, 
 a prime element of $P_F$, the residue field of $F$ and the valuation of $F$ 
 normalized by $v_F (\varpi_F) = 1$. 
  We set $q=q_F$ to be the number of elements in $k_F$. 
  Let $W_F$be the absolute Weil group of $F$.
 Hereafter we fix an additive character $\psi$ of $F$ whose conductor is
  $P_F$, i.e., $\psi$ is trivial
on $P_F$ and not trivial on ${\mathcal O}_F$. 
  For an extension $E$ over $F$, we denote by $\tr_E$, $\n_E$ the trace and norm
to $F$ respectively. We set $\psi_E=\psi\circ\tr_E$
and $\chi_E=\chi\circ\n_E$ for a quasi-character $\chi$ of $F^{\times}$.
Let $\theta$ be a quasi-character of $E^{\times}$. We denote by  $f(\theta)$ 
an integer such that $1+P_E^{n+1}\not \subset \Ker \theta$ and 
$1+P_E^n \subset \Ker\theta$.  
The Gauss sum $G(\theta,\psi_E)$ is defined by
\begin{equation}\label{gauss sum}
G(\theta,\psi_E)=
\begin{cases}
q_E^{-1/2}\sum_{x\in k_E^{\times}}\theta^{-1}(x)\psi_E(x) & \text{if $f(\theta)=1$} \\
q_E^{-1/2}\sum_{x\in k_E}\theta^{-1}(1+\varpi_E^m x)\psi_E(\varpi_E^m x) & 
\text{if $f(\theta)=2m+1$}.
\end{cases}
\end{equation}
The $\lambda$-factor $\lambda_{E}$ of $E/F$ is
defined by
\begin{equation}\label{lambda}
\lambda_{E}=\dfrac{\varepsilon(\Ind_{W_E}^{W_F} 1_{W_K},s,\psi_F)}
{\varepsilon( 1_{W_K},\psi_E)}.
\end{equation}
It is well-known that
\begin{equation}\label{lambda2}
\varepsilon(\Ind_{W_E}^{W_F} \sigma,s,\psi_F)=\lambda_{E}^{\dim \sigma}
\varepsilon(\sigma,s,\psi_E)
\end{equation}
for any representation $\sigma$ of $W_E$.
The trace of matrix is denoted by $\Tr$.
For an irreducible admissible representation $\pi$ of $\GL_l(F)$,        
the conductoral exponent of $\pi$ is defined to be the integer $f(\pi)$ such
that the local constant $\varepsilon (s, \pi ,\psi)$ of Godement-Jacquet 
\cite{GJ} is the form $a q^{-s(f(\pi ) - l)}$.

  Let $G$ be a totally disconnected, locally compact group. We denote 
 by $\widehat G$ the set of (equivalence classes 
 of) irreducible admissible representations of $G$.
  For a closed subgroup $H$ of $G$ and a representation $\rho$ of $H$, 
  we denote by $\Ind_H^G\rho$
 (resp. $\ind_H^G\rho$) the induced representation 
 (resp. compactly induced representation) of $\rho$ 
 to $G$.  For a representation $\pi$ of $G$, we denote by $\pi\vert_H$ 
 the restriction of $\pi$ to $H$.

\section{Construction of the representation $\mathrm{GL}_l(F)$}
Let $l$ be an arbitrary prime number (we allow the case $l=p$).
We set $V_F=F^l$ so that $\M_l(F)=\End_F(V_F)$
and $\GL_l(F)=\Aut_F(V_F)$. Throughout this paper, 
we write $G=G_{F}=\GL_l(F)$ and $G_{K}=\GL_l(K)$
In this section, we review the construction of the supercuspidal 
representation of $\GL_l(F)$ and its lift to $\GL_{l}(K)$ where $K/F$ is
a tamely ramified  extension. Most of the contents of this section are well-known
 (See \cite{Ca},\cite{Mo} and \cite{BH2}.)


\begin{defn}
Let $\calL=\{L_i\}_{i \in {\mathbb Z}}$ be the set of $\of$-lattices in $V_F$.
$\calL$ is said to be a uniform lattice chain of
period $e=e(\calL)$ if the following conditions hold for all $i\in{\mathbb Z}$:

\begin{enumerate}
\item
$L_{i+1} \subset L_i$,
\item
$P_F L_i = L_{i+e}$,
\item
$\dim_{k_F}(L_i/L_{i+1})=l/e$.
\end{enumerate}
\end{defn}

Since we assume $l$ is a prime, the period $e({\mathcal L})$ is either $l$ or $1$.

\begin{defn}
For a uniform lattice chain $\calL=\{L_i\}_{i \in {\mathbb Z}}$ of period
$e$, we set 
 $$
  \fA(\calL)=\{f \in \M_l(F) | f( L_i) \subset L_i \ \text{for all $i$}\},
 $$
 Then $\fA(\calL)$ is a principal order in $\M_l(F)$ and its Jacobson radical
 $\fP(\calL)$ is
 $$
 \{f \in \M_l(F) | f( L_i) \subset L_{i+1} \ \text{for all $i$}\}.
 $$
 We also set the period $e(\fA)$ of $\fA$ is the period of $\calL$.
 Put  $U(\fA)=\fA^{\times}$ ,$U(\fA)^n=1+\fP^n$ for any positive integer $n$ and
 $$
 \fK(\fA)=\Aut(\calL)=\{x \in \GL_l(F) |  x^{-1}\fA x =\fA \}.
$$
 \end{defn}

By taking an appropriate $\of$-basis of $L_0$,  we express the principal orders by the
following matrix form.
If  $e(\calL)=l$,  $\fA$ (resp. $\fP(\calL)$) is $G$-conjugate to 
$\M_l(\of)$ (resp. $\M_l(P_F)$). When $e(\calL)=1$,  up to $G$-conjugacy,
$$
\fA(\calL) =   \left\{
\begin{pmatrix}
\of & \of & \cdots & \of \\
P_F & \of & \cdots & \of \\
\hdotsfor{4} \\
P_F & P_F & \cdots & \of
\end{pmatrix}
\right\}
$$
and
$$
{\fP}(\calL) =   \left\{
\begin{pmatrix}
P_F & \of & \cdots & \of \\
P_F & P_F& \cdots & \of \\
\hdotsfor{4} \\
P_F & P_F & \cdots & P_F
\end{pmatrix}
\right\}.
$$

Let $r,n$ be integers satisfying
 $$
 n>r\geq \left[\frac{n}{2}\right]\geq 0,
 $$
where $[x]$ denote the greatest integer $\leq x$.  
For $\beta \in \M_l(F)$, we define a function $\psi_\beta$ on $U(\fA)^r$ by
\begin{equation}
\psi_\beta(1+x)=\psi(\Tr \beta x).
\end{equation}
Then the map $u \mapsto \psi_\beta$ induces an isomorphism between 
${\fP}^{-r+1}/{\fP}^{-n+1}$ and the complex dual, $(U(\fA)^r/U(\fA)^n)\sphat$,  of $U(\fA)^r/U(\fA)^n$.

\begin{defn}
Let $E/F$ be a field extension  in $\M_l(F)$.
 An element $\beta \in E$ is said to be $E/F$-minimal if the following conditions hold:

\begin{enumerate}
\item
$(v_E(\beta),e(E/F))=1$. 
\item
$k_F(\varpi_F^{-v_E(\beta)}\beta^{e(E/F)} \mod {P_E})=k_E$.
\end{enumerate}
\end{defn}

When $E\subset \M_L(F)$  and $E\neq F$, $E/F$ is an extension of degree $l$ since
$l$ is a prime. Thus we can identify $E$ with $V_F$. By this identification, 
$\{P_E^i\}_{i \in {\mathbb Z}}$ becomes a uniform lattice chain of period $e(E/F)$ . We put
$\fA(E)=\fA({P_E^i})$.

\begin{prop}\label{trivial intertwine}
Suppose $\beta$ is $E/F$-minimal  and $E\neq F$. For $\fA=\fA(E)$, we have:
\begin{enumerate}
\item
$\fK(\fA)=E^{\times}U(\fA)$ and $E^{\times} \cap U(\fA) = \Oe^{\times}$.
\item
$E \cap \fP^m =P_E^m$ for all integers $m$ and $E^{\times} \cap U(\fA)^m =
1+P_E^m$ for all integers $m \geq 1$.
\item
Let $x \in \fP^l$. If $\beta x-x\beta \in \fP^{m+l+1}$, then $x \in E+\fP^{l+1}$.
\end{enumerate}
\end{prop}
\begin{proof}
The last assertion of the above proposition is due to Carayol (see \cite{Ca}).
The rest is obvious.
\end{proof}

We shall construct the irreducible supercuspidal representations of $\GL_l(F)$ from $E/F$-minimal
elements. Let $E/F$ be a field extension of degree $l$, $\beta$ an $E/F$-minimal element and
 $\fA=\fA(E)$.
Put $v_E(\beta)=1-n<0$ and $m=[n/2]$. Then $\psi_\beta$ 
is a quasi-character of $U(\fA)^m$ whose kernel contains $U(\fA)^n$. 
 Put $H=E^{\times}U(\fA)^m$ and define 
a quasi-character $\rho_{\beta,\theta}$ of $H$ by 
\begin{equation}
   \rho_{\beta,\theta}(h \cdot g)=\theta(h)\psi_{\gamma}(g) 
    \qquad \text{for} \quad h \in E^{\times},\quad 
    g \in U(\fA)^m
\label{rho}
\end{equation}
where $\theta$ is a quasi-character of $E^{\times}$ such that 
$\theta|_{E^{\times}\cap U(\fA)^m}=\psi_\beta|_{E^{\times}\cap U(\fA)^m}$.
We note $f(\theta)=1-v_E(\beta)=n$ where $f(\theta)$ 
is the exponent of the conductor of $\theta$ i.e. the 
minimum integer such that $\Ker \theta\subset 1 + P_E^{n}$.

Put $J$ be the normalizer of $\psi_\beta$ in $\fK(\fA)$ i.e. 
$$   
J=\{a \in \fK(\fA) \ | \ \psi_\beta^a=\psi_\beta\}
$$
where $\psi_\beta^a(x)=\psi_\beta(a^{-1}xa)$ for $x \in U(\fA)^m$.
Then $J=E^{\times}U(\fA)^{m'}$ where $m'=[n/2]$ by virtue of Proposition~\ref{trivial intertwine}.
Put $\eta_{\beta,\theta}=\Ind_H^{\fK(\fA)} \rho_{\beta,\theta}$.

When $n$ is even, i.e. $n=2m$, then $J=H=E^{\times}U(\fA)^m$.
By the Clifford theory,
$\eta_{\beta,\theta}$ is an irreducible representation of $\fK(\fA)$.
We put
\begin{equation}
\kappa_{\beta,\theta}=\eta_{\beta,\theta}.
\label{kappa1} 
\end{equation}

When $n$ is odd, i.e. $n=2m-1$, then $J=E^{\times}U(\fA)^{m-1}$.
Thus $\eta_{\beta,\theta}$ is not an irreducible representation of $\fK(\fA)$.
Even  in this case, we can describe the irreducible component of $\eta_{\beta,\theta}$ by $\beta$ and $\theta$.
If $E/F$ is unramified. we put
\begin{equation}
\kappa_{\beta,\theta}=
\dfrac{(-1)^l(q^{l(l-1)/2}-(-1)^{l-1})(q-1)}{q^{l(l-1)/2}(q^l-1)}
\sum\limits_{\chi \in (E^{\times}/F^{\times}(1+P_E))\sphat}
 \eta_{\beta,\theta\otimes\chi}+(-1)^{l-1}
 \eta_{\beta,\theta}.
\label{kappa0}
\end{equation}
Now we assume we treat the case $E/F$ is ramified.
If  $l \neq p$, we put
\begin{equation}
\kappa_{\beta,\theta}=
\dfrac{1-\left(\dfrac{q}{l}\right)q^{(l-1)/2}}{lq^{(l-1)/2}}
\sum\limits_{\chi \in (E^{\times}/F^{\times}(1+P_E))\sphat}
 \eta_{\beta,\theta\otimes\chi}+\left(\dfrac{q}{l}\right)
 \eta_{\beta,\theta}
\label{kappa2}
\end{equation} 
where $\left(\dfrac{q}{l}\right)$ is the Legendre symbol. By Lemma~3.5.30 and Lemma~3.5.33 in \cite{Mo}, 
the virtual representation $\kappa_{\beta,\theta}$ is an irreducible component of $\eta_\theta$.

Next we treat the case $l=p$.  If $f$ is odd, we put
\begin{equation}
\kappa_\theta=
\sum\limits_{i=0}^{p-1} 
\left(\dfrac{1}{pq^{(l-1)/2}}+\left(\dfrac{i}{p}\right)\right)\dfrac{p^{(f-1)/2}}{G_0G(\beta)}
\eta_{\theta\otimes\chi^i}
\label{kappa3}
\end{equation} 
where $\chi$ is a generator of $(E^{\times}/F^{\times}(1+P_E))\sphat$ determined by 
$\chi(\varpi_E)=\exp(2\pi\sqrt{-1}/p)$ and 
$G_0,G(\beta)$ are  Gauss sums defined by 
\begin{align}
G(\beta) &= \frac{1}{\sqrt{q}}\sum_{x \in k_E}\psi(\tr_{k_E/k_F}\frac 12 
\beta\varpi_E^{2(m-1)}(-1)^{(p+1)/2}x^2) \\
G_0&=\frac{1}{\sqrt{p}}\sum_{a=1}^{l}\left(\frac ap\right)\exp(2\pi\sqrt{-1}a/p).
\label{gauss sums}
\end{align}
When $f$ is even, we put
\begin{equation}
\kappa_\theta=
\sum\limits_{\chi \in (E^{\times}/F^{\times}(1+P_E))\sphat}
\dfrac{q^{1/2}G(\beta)-q^{(p-1)/2}}{G(\beta)pq^{p/2}}\eta_{\theta\otimes\chi}+\dfrac{1}
{q^{1/2}G(\beta)}\eta_{\chi}.
\label{kappa4}
\end{equation} 
By Proposition~3.4 in \cite{T}, $\kappa_\theta$ is an irreducible component of $\eta_\theta$.

Finally we consider the level $1$ supercuspidal representation.
Let $E/F$ be an unramified extension of degree $l$ ,$\theta$  a quasi-character of $E^{\times}$ which is
trivial on $1+P_E$ and $\fA=\fA(E)$. Then there is an irreducible representation 
 $\kappa_\theta$
 of $U(\fA)$ which is trivial on $U(\fA)^1$ and its tensor product with the
 pull-back of the Steinberg representation of 
 $U(\fA)/U(\fA)^1\simeq \GL_l(k_F)$ is the representation induced
 by the one-dimensional representation 
 $tx \mapsto\theta(t),t\in \Oe^{\times},x\in U(\fA)^1$, of the subgroup 
 $\Oe^{\times}U(\fA)^1$. We denote by $\kappa_\theta$ the representation 
 $tx\mapsto\theta(t)\kappa_\theta(x),t \in F^{\times}, x\in U(\fA)$,
 of $\fK(\fA)$.

\begin{thm}
Let the notation be as above.
Then $\kappa_{\beta,\theta}$ and $\kappa_\theta$ are  irreducible representations of $\fK(\fA)$. Put $\pi_F(\beta,\theta)=\ind_{\fK(\fA)}^G \kappa_{\beta,\theta}$ and 
$\pi_F(\theta)=\ind_{\fK(\fA)}^G \kappa_\theta$.
Then $\pi_F(\beta,\theta)$  and $\pi_F(\theta)$ are irreducible
supercuspidal representations of $G=\GL_{l}(F)$ with $f(\pi_F(u,\theta)=f(E/F)(f(\theta)-1)+l$ and $f(\pi_F(\theta))=l f(\theta)$. Every irreducible supercuspidal representation of $G$ can be written in
the form $\chi\pi_F(\beta,\theta)$ or $\chi\pi_F(\theta)$ for some $E/F$-minimal element $\beta$, $\theta\in \widehat{F(\beta)^{\times}}$ and 
$\chi\in \widehat{F^{\times}}$.
\end{thm}

The $\varepsilon$-factors of all supercuspidal representations of $G$ have been calculated completely.
(See \cite{Mo},\cite{KM}).
\begin{thm}\label{epsilon}
Let $\pi_F(\beta,\theta)$ and $\pi_F(\theta)$ be as above. Put $n=f(\theta)$.
For $\chi \in \widehat{F^{\times}}$,  we pick an element $c_\chi \in F$ such
that $\chi(1+x)=\psi_F(c_\chi x)$ for $x\in P_F^{[(f(\chi)+1)/2]}$. (If $f(\chi)\leq 1$, we take $c_\chi=0$). 
Put  $n(\chi)=\max (n,e(E/F)(f(\chi)-1)+1)$ and $\beta_\chi=\beta+c_\chi$. 
\begin{enumerate}
\item
 If $n(\chi)$ is even, 
$$
\varepsilon(\chi\pi_F(\beta,\theta) ,s,\psi)=\psi_E(\beta_\chi)(\chi_E\theta)(\beta_\chi^{-1})|\beta_\chi|_E^{s}.
$$

\item
If $n(\chi)=n=1$, 
$$
\varepsilon(\chi\pi_F(\theta) ,s,\psi)=(-1)^{l-1}\varepsilon(\chi_E\theta,s,\psi_E).
$$

\item
If $n(\chi)\neq 1$ is odd,
$$
\varepsilon(\chi\pi_F(\beta,\theta),s,\psi)=
\psi_E(\beta_\chi)(\chi_E\theta)(\beta_\chi^{-1})|\beta_\chi|_E^{s}G
$$
where the Gauss sum $G$ is defined by
$$
G=
\begin{cases}
G(\theta,\psi_E) & \text{if} \quad n=n(\chi) \text{and $E/F$is tamely ramified}\\
G(\beta) & \text{if} \quad n=n(\chi) \text{and $E/F$ is wildly ramified} \\
\lambda_{E}G(\chi,\psi_F)^l & \text{if} \quad n>n(\chi)
\end{cases}
$$
where $\lambda_{E}$ is defined in \eqref{lambda}.
\end{enumerate}
\end{thm}

\section{Explicit correspondences and tame base change lift} 
Now we consider some correspondences between
 ${\mathcal A}_F(l)$ and ${\mathcal G}_F(l)$ which satisfy the conditions (i)-(iv) of
the local Langlands correspondence. When $l\neq p$ or $l=p$ and $E/F$ is unramified, this is a special case of Howe-Moy correspondence.

\begin{defn}
A quasi-character $\theta$ of $E^{\times}$ is called generic if
 $f(\theta) \not\equiv 1 \bmod{l}$. For a generic character $\theta$ of
$E^{\times}$,   $\beta_\theta 
  \in P_E^{1-f(\theta)}-P_E^{2-f(\theta)}$ is defined by
  \begin{equation}
   \theta (1+x) = \psi_E (\beta_\theta x) \quad \text{for} \quad x 
   \in P_E^{[(f(\theta)+1)/2]}.
\label{gamma}  
\end{equation}
Then $\beta$ is $E/F$-minimal.
We denote by $\widehat{E_{gen}^{\times}}$ the set of generic quasi-characters
of $E^{\times}$.
\end{defn}

\begin{rem}\label{tamepi}
When $E/F$ is tamely ramified, the generic quasi-character $\theta$ determines uniquely 
$\pi_F(\beta_\theta,\theta)$. (See \cite{Mo}).  In this case we simply denote $\pi_F(\beta,\theta)$ by $\pi_F(\theta)$. 
When $l=p$, we need $\beta$ to determine the representation $\pi_F(\beta,\theta)$.
\end{rem}

To separate the wildly ramified case, 
we introduce some notations. Let ${\mathcal A}_l^{wr}(F)$  denote the
set $\pi=\chi\pi_F(\beta,\theta)\in {\mathcal A}_l(F)$  
with the property that $F(\beta)/F$ is wildly ramified.
$\pi \in {\mathcal A}_l^{wr}(F)$ is equivalent to $l=p$ and $\pi\simeq \chi\pi$
for some  unramified quasi-character $\chi\neq 1$ of $F^{\times}$. 
We put ${\mathcal A}_l^{t}(F)=
{\mathcal A}_l(F)\setminus{\mathcal A}_l^{wr}(F)$. 
Similarly let ${\mathcal G}_l^{wr}(F)$  denote the
set $\sigma\in {\mathcal G}_l(F)$  with the property that $\sigma\otimes\chi$ is 
equivalent to $\pi$ for some unramified quasi-character $\chi\neq 1$ of $F^{\times}$ and $l=p$.
We also put  ${\mathcal G}_l^{t}(F)={\mathcal G}_l(F)\setminus{\mathcal G}_l^{wr}(F)$.
If $p\neq l$, ${\mathcal A}_l(F)={\mathcal A}_l^{t}(F)$ and 
${\mathcal G}_l(F)={\mathcal G}_l^{t}(F)$
The Howe-Moy correspondence gives a bijection between 
${\mathcal G}_l^{t}(F)={\mathcal A}_l^{t}(F)$. (See \cite{Mo} and \cite{Ge}.)

 
 If $E/F$ is tamely ramified, $\lambda_{E}$ is easily calculated.
\begin{lem}\label{lambda's value}
Let $E/F$ is a tamely ramified extension of degree $l$.
Then
$$
\lambda_{E}=
\begin{cases}
(-1)^{l-1} & \text{if $e(E/F)=1$}, \\
\left(\dfrac{q}{l}\right) & \text{if $e(E/F)=l\neq 2$} \\
q^{-1/2}\sum_{x \in k_E} \sgn_{E/F}^{-1}(x)\psi_E(x) &
\text{if $e(E/F)=l=2$} 
\end{cases}
$$
\end{lem}
\begin{proof}
See (2.5.3), (2.5.5) and Proposition 2.5.11 of \cite{Mo}.
\end{proof}

\begin{thm}\label{tameLLC}
Let $E$ be a tamely ramified extension of $F$ of degree $l$ 
and $\theta$ be a generic quasi-character of $E^{\times}$. 
We define a quasi-character $\delta_{E}$ of $E^{\times}$ as follows:

When $e(E/F)\neq 2$, $\delta_{E}(x)=\lambda_E^{v_E(x)}$.

When $e(E/F)=2$,  
$$
\delta_{E}(x)=
\begin{cases}
1 & {\text if}\quad  x\in1+P_E, \\
\sgn_{E/F}(x) & \text{if}\quad  x\in F^{\times}, \\
\lambda_E & \text{if}\quad  x=\beta_\theta.
\end{cases}
$$
We set 
$$
\sigma_F(\theta)= 
\delta_E \Ind_{W_E}^{W_F}.
$$
\begin{enumerate}
\item
the representation $\sigma_F(\theta)$ belongs to ${\mathcal G}_l^{t}(F)$  
and any element of ${\mathcal G}_l^{t}(F)$
 can be written in the form $\chi\sigma_F(\theta)$ for an extension $E/F$ of degree $l$, 
a generic character of $E^{\times}$ and a quasi-character $\chi$ of $F^{\times}$.
\item
Define the map $\Phi_l^F$ from ${\mathcal G}_l^{t}(F)$ to ${\mathcal A}_l^{t}(F)$ by
$$
\Phi_l^F(\chi\sigma_F(\theta)) =\chi\pi_F(\delta_{E}\theta).
$$
Then $\Phi_l^F$ is a bijection which satisfies the following conditions:
\begin{enumerate}
  \item 
  For $\chi \in \widehat{F^{\times}}$and $\sigma\in{\mathcal G}_l^{t}(F)$, 
  $$
  \Phi_l^F(\chi\sigma)=\chi\Phi_l^F(\sigma).
  $$
  \item 
  For $\sigma\in{\mathcal G}_l^{t}(F)$, 
  $$
  \Phi_l^F(\check{\sigma})=\Phi_l^F(\sigma)\spcheck .
  $$
  \item
  Let $\omega_\pi$ denote the central quasi-character of $\pi\in{\mathcal A}_l(F)$.  
  For $\sigma\in{\mathcal G}_l^{t}(F)$.
  $$
  \omega_{\Phi_l^F(\sigma)}=\det \sigma.
  $$
  \item
  For $\sigma\in{\mathcal G}_l^{t}(F)$, 
  $$
  \varepsilon(\Phi_l^F(\sigma),s,\psi_F)=\varepsilon(\sigma,s,\psi_F).
  $$
\end{enumerate}
\end{enumerate}
\end{thm}

Since ${\mathcal G}_p^{wr}(F)$ contains non-monomial representations, the correspondence between ${\mathcal G}_p^{wr}(F)$ and ${\mathcal A}_p^{wr}(F)$ becomes more complicated.
We use the tame lifting of Bushnell-Henniart \cite{BH2}. For any tamely ramified extension $K/F$,  including the case $K/F$ is non-Galois, the tame lifting map ${\mathbf l}_{K}$ from
${\mathcal A}_{p^i}^{wr}(F)$ to ${\mathcal A}_{p^i}^{wr}(K)$ is constructed by Bushnell-Henniart.  Since we consider the case $i=1$, this base change is easy to describe.
Since $K/F$ is tamely ramified, $E\otimes_F K=EK$ is an extension of field of $K$,
$G_K=G(K)$ can be identified with $\Aut_K(E\otimes_F K)$ and $\beta=\beta\otimes 1$ becomes
 an $EK/K$-minimal element in $V_K=\End_K(EK)$. 
 Moreover  if $\theta$ is a quasi-character of $E^{\times}$ such that $\theta(1+x)=\psi(\tr_{E/F}\beta x)$ for $x\in P_{E}^m$, then $\theta\circ n_{EK/E}(1+x)=\psi_{K}(\tr_{EK/K}\beta x)$ for $x\in P_{EK}^m$. Therefore
  we get an irreducible supercuspidal representation 
 $\pi_{K}(\beta,\theta\circ \n_{EK/E})\in {\mathcal A}_p^{wr}(K)$.
 
\begin{thm}\label{explicit lift}
Let $K/F$ be an extension of degree prime to $p$ and ${\mathbf l}_{K}$  the lifting from 
${\mathcal A}_{p}^{wr}(F)$ to ${\mathcal A}_{p}^{wr}(K)$ defined by (5.3.3) in \cite{BH2}. 
Put $\Delta_{K}=\det\Ind_{W_K}^{W_F} 1_{W_K}\in \widehat{F^{\times}}$ and $\tilde{\Delta}=\Delta_{K}\circ\n_{E/F} \in\widehat{E^{\times}}$.
For $\pi_F(\beta,\theta)\in {\mathcal A}_{p}^{wr}(F)$ and 
$\chi\in \widehat{F^{\times}}$, we have:
\begin{align*}
{\mathbf l}_{K}(\chi\pi_F(\beta,\theta)) &=
\chi_K\pi_K(\beta,(\tilde{\Delta}^{e(E/F)-1}\theta)\circ\n_{EK/E}) \\
&=
\begin{cases}
\chi_K\pi_K(\beta,\theta\circ\n_{EK/E}) & e(E/F)\neq 2\\
\chi_K\pi_K(\beta,(\tilde{\Delta}\theta)\circ\n_{EK/E}) & e(E/F)=2.
\end{cases}
\end{align*}
\end{thm}
\begin{proof}
Since two lifting maps are compatible with twist of quasi-character of $F^{\times}$,
we may assume $\chi=1$.
By Proposition 10.2 of \cite{BH2}, it suffices to say
$$
\varepsilon({\mathbf l}_{K}(\pi_F(\beta,\theta)),s,\psi_K)=
\varepsilon(\pi_K(\beta,(\tilde{\Delta}^{e(E/F)-1}\theta)\circ\n_{EK/E}),s,\psi_K).
$$
(Other conditions (a) and (b) in Proposition 10.2 of \cite{BH2} are obvious in our case.)
Theorem 1.6 of \cite{BH2}  tells us that
$$
\lambda_K^p\varepsilon({\mathbf l}_{K}(\pi_F(\beta,\theta)),s,\psi_K)=
\Delta(\n_{E/F}(\beta))\varepsilon(\pi_F(\beta,\theta),s,\psi_F)^{[K:F]}.
$$
On the other hand, it follows from Proposition 2.2.11 of \cite{KM} that 
$$
\lambda_K^p\varepsilon(\pi_K(\beta,\theta\circ\n_{EK/E}),s,\psi_K)=
\Delta(\n_{E/F}(\beta))\varepsilon(\pi_F(\beta,\theta),s,\psi_F)^{[K:F]}
$$
if $p \neq 2$. 
(Proposition 2.2.11 of \cite{KM} assumes $K/F$ is Galois, but it holds including the
case $K/F$ is non-Galois since Proposition 2.5.16 of \cite{Mo}  holds for any
tamely ramified extension $K/F$.)  Hence the assertion holds when $p\neq 2$.
When $p=2$, $n(\pi_F(\beta,\theta)=1-v_E(\beta)$ is even since $(v_E(\beta),p)=1$.
Therefore Theorem~\ref{epsilon} tells us
$$
\varepsilon(\pi_K(\beta,\theta\circ\n_{EK/E}),s,\psi_K)=
\varepsilon(\pi_F(\beta,\theta),s,\psi_F)^{[K:F]}.
$$
Since $\Delta_{K}\circ \n_{K/F}$ is unramified and $\Delta_{K}^{-1}=\Delta_{K}$,
$$
\varepsilon(\pi_K(\beta,(\tilde{\Delta}\theta)\circ\n_{EK/E}),s,\psi_K)=
\Delta_{K}\circ\n_{K/F}(\beta)\varepsilon(\pi_K(\beta,\theta\circ\n_{EK/E}),s,\psi_K).
$$
Hence our assertion holds.
\end{proof}

\begin{rem}\label{delta and Delta}
Two quasi-characters $\Delta_{K}$ and  $\delta_{K}$ is closely related.
If $e(K/F)$ is odd.,  $\Delta_{K}\circ\n_{K/F}=\delta_{K}$. (See Corollary~2.5.5 of \cite{Mo}.)
\end{rem}

Using the tame lifting map ${\mathbf l}_{K}$, Bushnell-Henniart (\cite{BH2}) has constructed
the correspondence ${\mathcal G}_{p^i}^{wr}(F)$ to ${\mathcal A}_{p^i}^{wr}(F)$.
For $i=1$, this map coincides with the local Langlands correspondence and
is compatible with ${\mathbf l}_{K}$. This follows as a special case of Lemma 5.2 in \cite{BH3}.

\begin{prop}\label{compatibility of lifting}
Let $\Lambda_l^F$ be the local Langlands map. 
Then for any tamely ramified extension $K/F$ and 
$\sigma\in{\mathcal G}_{p}^{wr}(F)$, we have:
$$
{\mathbf l}_{K}\Lambda_l^F(\sigma)=\Lambda_l^K(\sigma |_{W_K}).
$$
\end{prop}

\begin{proof}
By Lemma 5.2 in \cite{BH3}, it suffices to say that the exponent $f(\pi_{\beta,\theta})$ 
of the conductor of $\pi_{\beta,\theta} \in {\mathcal A}_{p}^{wr}(F)$ is prime to $p$. It follows from the fact that $f(\pi_{\beta,\theta} )\equiv -v_E(\beta) \mod{p}$.
\end{proof}

We define the lift ${\mathbf l}_{K}$ for $\pi\in {\mathcal A}_l^{t}(F)$ as in the case
$\pi\in {\mathcal A}_l^{wr}(F)$.

\begin{defn}
Let $E/F$ be  an extension of degree $l$ , $\theta\in\widehat{E^{\times}_{gen}}$
and $\chi\in\widehat{F^{\times}}$. Assume $K$ is a tamely ramified extension of $F$ such that
$([K:F],l)=1$. Then we define ${\mathbf l}_{K}(\chi\pi_F(\theta))$ by 
\begin{align*}
{\mathbf l}_{K}(\chi\pi_F(\beta,\theta)) & =
\chi_K(\Delta_{K}\circ \n_{K/F})^{e(E/F)-1}\pi_K(\beta,\theta\circ\n_{EK/E}) \\
&=
\begin{cases}
\chi_K\pi_K(\beta,\theta\circ\n_{EK/E}) & e(E/F)\neq 2\\
\chi_K(\Delta_{K}\circ \n_{K/F})\pi_K(\beta,\theta\circ\n_{EK/E}) & e(E/F)=2.
\end{cases}
\end{align*}
\end{defn}

This lifting is compatible with $\Phi_l$.

\begin{prop}\label{compatibility of lifting II}
Let $K/F$ be a finite, tamely ramified extension satisfying $K \cap E=F$.
For $\sigma\in{\mathcal G}_{l}^{t}(F)$,
$$
{\mathbf l}_{K}\Phi_l^F(\sigma)=\Phi_l^K(\sigma |_{W_K}).
$$
\end{prop}

\begin{proof}
Since ${\mathbf l}_{K}$ and  $\Phi_l$ are compatible with quasi-character twist,
we may assume $\sigma=\sigma_F(\theta)$ for $\theta\in \widehat{E^{\times}_{gen}}$.
By the definition of ${\mathbf l}_{K}$ and  $\Phi_l$,
$$
(\Phi_l^K)^{-1}({\mathbf l}_{K}(\Lambda_l^F(\sigma_F(\theta))))=
\Ind_{W_{EK}}^{W_{K}}\delta_{EK/K}((\tilde{\Delta}^{e(E/F)-1}\theta)\circ\n_{EK/K}).
$$
On the other hand, it follows from $W_EW_K=W_F$ and $W_E \cap W_K=W_{EK}$ that
Mackey's Theorem tells us
$$
\sigma_F(\theta)|_{W_K}=\Ind_{W_{EK}}^{W_{K}}(\delta_{E}\theta)\circ \n_{EK/E}.
$$
Thus it suffices to say that 
\begin{equation}\label{delta's equation}
\delta_{EK/K}(\tilde{\Delta}^{e(E/F)-1}\circ\n_{EK/K})=
\delta_{E}\circ\n_{EK/E}.
\end{equation}
When $e(E/F)$ is odd, this is obvious. So we assume $e(E/F)=2$.
By the definition of $\delta$, we have only to show \eqref{delta's equation} for $x \in K^{\times}$
and $\beta$.
For $x\in K^{\times}$, 
\begin{align*}
\delta_{EK/K}(\tilde{\Delta}\circ\n_{EK/K}(x))&=
\delta_{EK/K}(x)\Delta_{K}\circ\n_{K/F}(x^2) \\
&=\sgn_{EK/K}(x).
\end{align*}
since $\Delta_{K}$ has at most order $2$. The right hand side of \eqref{delta's equation} 
becomes $\sgn_{E/F}(\n_{EK/E}(x))$, which equals to $\sgn_{EK/K}(x)$ since $[K:F]$ is odd.
We compare the value of both sides of \eqref{delta's equation} at $\beta$.
The left hand side is $\delta_{EK/K}(\beta)\Delta_{K}(\n_{E/F}(\beta))^{[K:F]}$.
It follows from  Remark~\ref{delta and Delta} that $\Delta_{K}(\n_{E/F}(\beta))=
\delta_{K}^{v_K(\beta)}$. Thus it amounts to $\lambda_{EK/K}\lambda_{K}$.
The right hand side becomes $\lambda_{E}^{[K:F]}$. After all, the equation
$\lambda_{EK/E}\lambda_{K}=1$ gives the result. 
When $[K:F]$ is prime, it follows from Lemma~\ref{lambda's value}.  
The composite case is obtained by the transitivity property 
of $\lambda$-factor.
\end{proof}

We need to show that the Howe-Moy correspondence $\Phi_{l'}$ coincides with 
the Local Langlands correspondence $\Lambda_{l'}$.

\begin{thm}\label{Phi=Lambda}
For any prime $l'\neq p$,
$$
\Phi_{l'}^F=\Lambda_{l'}^F.
$$
\end{thm}
\begin{proof}
If $l'=2$, it follows from Converse Theorem (\cite{CPS}). 
We assume $l'$ is an odd prime.
Let $\pi \in {\mathcal A}_F(l')$ . Then there exist
an extension  $E/F$ of degree $l'$, $\theta\in \widehat{(E^{\times})_{gen}}$ and $\chi\in\widehat{F^{\times}}$ such that $\pi=\chi_{E}\pi_F(\theta)$ as in Remark~\ref{tamepi}.
When $E/F$ is unramified, Theorem 9.2 (\cite{Spi}) implies $\Phi_{l'}^F(\pi)=\chi_E\Ind_{W_E}^{W_F}\theta=\Lambda_{l'}^F(\pi)$. 
When $E/F$ is ramified, the assertion follows from Theorem B in \cite{BH6}.
\end{proof}
\begin{rem}
Theorem 9.2 (\cite{Spi}) is proved under the assumption $p>l$, but this assumption is dispensable.
The key point  is to prove that
$$
\Theta_{\pi}^{\kappa}(x)=\Theta_{\pi}(x) \quad \text{for $x \in E^{\times}\backslash F^{\times}(1+P_E^r)$}
$$
where $\Theta_\pi$ is a distribution character of $\pi$ and $\Theta_\pi^\kappa$ is a
$\kappa$-twisted distribution character of $\pi$ for $\kappa \in (F^{\times}/n_{E/F}(E^{\times}))
\sphat$. This is proved in Theorem~6.1 (\cite{Spi}) without using the assumption $p>l$.
\end{rem}
By Propostion~\ref{compatibility of lifting}, Propostion~\ref{compatibility of lifting II} and
Theorem~\ref{Phi=Lambda}, $\Phi_l$ is comatible with ${\mathbf l}$ for any prime $l$.
\begin{cor}\label{comatibility of lifting III}
Let $K/F$ be a finite, tamely ramified extension satisfying $K \cap E=F$.
For any prime $l$ and $\sigma\in{\mathcal G}_{l}(F)$,
$$
{\mathbf l}_{K}\Phi_l^F(\sigma)=\Phi_l^K(\sigma |_{W_K}).
$$
\end{cor}
\section{$\varepsilon$-factor of pairs}
In this section, we consider the $\varepsilon$-factor 
$\varepsilon(\pi_1\times\pi_2,s,\psi_F)$.
Let $l'$ be a prime not equal to $l$ and $p$. 
We treat the case $\pi_1 \in {\mathcal A}_F(l)$ and $\pi_2 \in {\mathcal A}_F(l')$.
Since the local Langlands correspondence and the Bushnell-Henniart base change lift are compatible with quasi-character twists, we may assume $\pi_1$ and $\pi_2$ are minimal.

\begin{thm}\label{main}
Let $\pi_1 \in {\mathcal A}_F(l)$ and $\pi_2 \in {\mathcal A}_F(l')$ where $l'$ is a prime
not equal to $l$ and $p$. 
Let $E_2/F$ be an extension of degree $l'$,  $\theta_2\in \widehat{(E_2^{\times})_{gen}}$ and $\chi_2\in\widehat{F^{\times}}$ such that $\pi_2=(\chi_2)_{E_2}\pi_F(\theta_2)$ as in Remark~\ref{tamepi}.
Then we have
\begin{equation}\label{Frob recip}
\varepsilon(\pi_1\times\pi_2,s,\psi_F)=
\lambda_{E_2}\varepsilon(\chi_2\delta_{E_2}\theta_2
{\mathbf l}_{E_2}(\pi_1),s,\psi_{E_2}).
\end{equation}
\end{thm}

\begin{proof}
It follows from $\Phi_{l'}^F=\Lambda_{l'}^F$ that 
$$
(\Lambda_{l'}^F)^{-1}(\pi_F((\chi_2)_{E_2}(\theta_2))=
\Ind_{W_{E_2}}^{W_F}(\chi_2\delta_{E_2}\theta_2).
$$
Put $(\Lambda_l^F)^{-1}(\pi_1)=\sigma_1$. Then we have:
\begin{equation*}
\varepsilon(\pi_1\times \pi_2,s,\psi_F)=
\varepsilon(\sigma_1\otimes
\Ind_{W_{E_2}}^{W_F}(\chi_2\delta_{E_2}\theta_2),s,\psi_F). 
\end{equation*}
Since
$$
\Ind_{W_{E_2}}^{W_F}\sigma_1\otimes\chi_2\delta_{E_2}\theta_2
=\Ind_{W_{E_2}}^{W_F}(\sigma_1|_{W_{E_2}}\otimes \chi_2\delta_{E_2}\theta_2)
$$
and 
$$
\varepsilon(\Ind_{W_{E_2}}^{W_F}\sigma,s,\psi_F)=
\lambda_{E_2}^{\dim \sigma}\varepsilon(\sigma,s,\psi_{E_1}) \quad \text{for}\quad \sigma\in{\mathcal G}_{E_2}(l'),
$$
we obtain
\begin{align*}
\varepsilon(\pi_1\times \pi_2,s,\psi_F)&=\varepsilon({\sigma_1}|_{W_{E_2}}\otimes
\Ind_{W_{E_2}}^{W_F}(\chi_2\delta_{E_2},s,\psi_F) \\
&=\lambda_{E_2}\varepsilon({\sigma_1}|_{W_{E_2}}\otimes
(\chi_2\delta_{E_2}\theta_2),s,\psi_{E_1}).
\end{align*}
Assume $l\neq p$, then $\pi_1$ can be written in the corm $\chi_1\pi_{\theta_1}$and 
$\sigma_2=\Ind_{W_{E_1}}^{W_F}(\chi_1\delta_{E_1}\theta_1$. 
By the Mackey decomposition and $W_{E_1}W_{E_2}=W_F$, we have
$$
(\Ind_{W_{E_1}}^{W_F}\chi_1\delta_{E_1}\theta_1)|_{W_{E_2}}=
\Ind_{W_{E_1E_2}}^{W_{E_2}}(\chi_1\delta_{E_1}\theta_1)\circ\n_{E_1E_2/E_1}.
$$ 
Since $(\chi_1\delta_{E_1}\theta_1)\circ\n_{E_1E_2/F}$ does not factor through 
$\n_{E_1E_2/E_2}$, 
$$
\Ind_{W_{E_1E_2}}^{W_{E_2}}((\chi_1\delta_{E_1}\theta_1)\circ\n_{E_1E_2/E_1} \in {\mathcal G}_{E_2}(l).
$$
Therefore we have :
\begin{align*}
\varepsilon(\pi_1\times\pi_2,s,\psi_F)&=\lambda_{E_2}
\varepsilon(\Ind_{W_{E_1E_2}}^{W_{E_2}}(\chi_1\delta_{E_1}\theta_1)\otimes
\chi_2\delta_{E_2}\theta_2\circ\n_{E_1E_2/E_2},s,\psi_{E_1}) \\
&=\lambda_{E_2}
\varepsilon(\chi_2\delta_{E_2}\theta_2\otimes(\chi_1)_{E_1E_2}
\pi_{E_2}((\theta_1\circ\n_{E_1E_2/E_1}),s,\psi_{E_2}).
\end{align*}
(The last equality follows from $\Lambda_{E_2}(\pi_{\theta_1}\circ\n_{E_1E_2/E_1})
=\Ind_{W_{E_1E_2}}^{W_{E_2}}(\theta_1\delta_{E_1})\circ\n_{E_1E_2/E_1}$.)

When $l=p$,  it follows from Proposition~\ref{compatibility of lifting} and 
Proposition~\ref{explicit lift} that
\begin{align*}
\Lambda_{E_2}({\sigma_1}|_{W_{E_2}}) &={\mathbf l}_{E_2}(\pi_1) \\
&=(\chi_1)_{E_1E_2}\pi_{E_2}(\beta_1,\theta_1\circ\n_{E_1E_2/E_1}).
\end{align*}
Thus we have
$$
\varepsilon(\pi_1\times \pi_2,s,\psi_F)=
\lambda_{E_2}\varepsilon(\chi_2\delta_{E_2}\theta_2\otimes (\chi_1)_{E_1E_2}
\pi_{E_2}(\beta_1,\theta_1\circ\n_{E_1E_2/E_1}),s,\psi_{E_2}).
$$
\end{proof}

By combining Theorem~\ref{epsilon} and Theorem~\ref{main}, 
we get the complete list of $\\varepsilon(\pi_1\times\pi_2,s,\psi_F)$ for $
\pi_1 \in {\mathcal A}_F(l)$ and $\pi_2\in {\mathcal A}_F(l')$ where $l$ is any prime and
$l'$ is  a prime $\neq l$.
 
\begin{rem}
By the result of \cite{BH6} , Theorem~\ref{main} may be extended to the case 
$\pi_1 \in {\mathcal A}_F^{t}(m)$ and $\pi_2\in {\mathcal A}_F^{t}(n)$ where
$ (m,n)=1$.  
\end{rem}


\end{document}